\documentclass{article}
\usepackage{graphics} 
\usepackage{graphicx}
\usepackage{epsfig} 
\usepackage{amsmath} 
\usepackage{amssymb}  
\usepackage{color}  

\usepackage{bm}
\usepackage{txfonts}

\newtheorem{proposition}{\textbf{Proposition}}
\newtheorem{lemma}{\textbf{Lemma}}
\newtheorem{theorem}{\textbf{Theorem}}
\newtheorem{remark}{\textbf{Remark}}

\newtheorem{definition}{\textbf{Definition}}
\newtheorem{conjecture}{\textbf{Conjecture}}
 \usepackage{comment}
\usepackage{subfigure}

\begin{document}
\title{\LARGE On the Meeting Time for Two Random Walks on a Regular Graph}

\author{Yizhen Zhang$^{\dag}$, Zihan Tan$^{\dag}$ and Bhaskar Krishnamachari$^{\ddag}$ \\ 
$^{\dag}$ Tsinghua University, Beijing, China\\ $^{\ddag}$ University of Southern California, Los Angeles, USA\\ \{yz-zhang11, zh-tan11 \}@mails.tsinghua.edu.cn, bkrishna@usc.edu}

\maketitle

\begin{abstract}
We provide an analysis of the expected meeting time of two independent random walks on a regular graph. For 1-D circle and 2-D torus graphs, we show that the expected meeting time can be expressed as the sum of the inverse of non-zero eigenvalues of a suitably defined Laplacian matrix. We also conjecture based on empirical evidence that this result holds more generally for simple random walks on arbitrary regular graphs. Further, we show that the expected meeting time for the 1-D circle of size $N$ is $\Theta(N^2)$, and for a 2-D $N \times N$ torus it is $\Theta(N^2 log N)$. 
\end{abstract}

\section{Introduction}


Consider a system of discrete-time random walks on a graph $G(V,E)$ with two walkers. Each time, they each independently move to a nearby vertex or stay still with given probabilities. Denote the transition matrix of a single walker by $P$, where $P(i,j)$ is the probability that one walker moves from $v_i$ to $v_j$ in a time slot. This process is assumed to start at steady state (i.e. uniform distribution) for each walker, and terminates when they meet at the same vertex. We denote this meeting time by $\tau$, which is a random variable with the expectation \textbf{E}$[\tau]$. Our objective is to analyze this quantity on $d$-regular graphs. 

\begin{figure}[hbp]
\centering
  \includegraphics[width=0.2\textwidth]{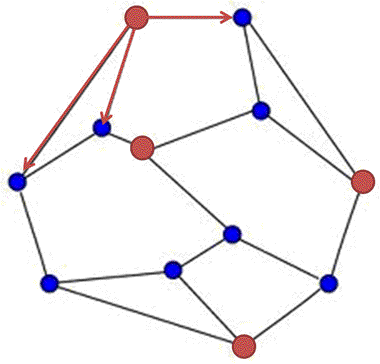}
  \caption{4 walkers on a 3-regular graph}
  \label{fig 1}
\end{figure}

It is instructive to consider the problem on the one-dimensional circle first. We study a circle with $N$ nodes, denoted by $V=\{0,1,2, \cdots ,N-1 \}$. The two walkers start from arbitrary position according to the initial distribution. Every step, the walker on $i$ chooses to move to $\{i-1,i+1\}$ (for simplicity of notation, assume that if $i=N-1$, then $i+1 = 0$ and similarly that if $i=0$ then $i-1 = N-1$ ) or stay still at $i$ with probability $\{p_1,p_2,p_3\}$ respectively.  

\begin{figure}[hbp]
\centering
  \includegraphics[width=0.6\textwidth]{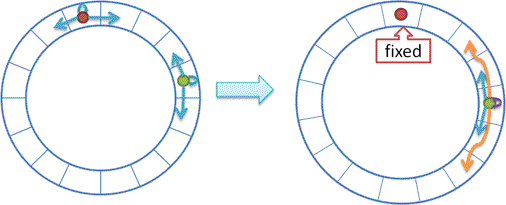}
  \caption{1-D circle}
  \label{fig 2}
\end{figure}

Since we are only concerned about the meeting time, the relative position of the two walkers is enough to describe that random variable. So we fix one walker at '$0$'. Then in this new equivalent model, the transition matrix of the other walker before the encounter is
$M = P \cdot P^T$.

A similar equivalent model can be defined for a $N \times N$ torus. Let $V = \{(x,y)|x,y=1,2,...,l\}$. Every step, the walker on $(x,y)$ moves to $(x\pm 1,y\pm 1)$ or stay still at $(x,y)$ with given probability. Define the index of $(x,y)$ to be \textbf{Ind}$(x,y)=(x-1)N+y$, we can get a $N^2$-order matrix $P$. Let $i,j$ denote the indices of two vertex $(x_i,y_i),(x_j,y_j)$. Then $P(i,j)$ denotes the probability that one walker moves from $(x_i,y_i)$ to $(x_j,y_j)$ each step. $P$ is a ``block-circulant matrix" defined in 3.1.2 .

Similar to the 1-D case, we fix one walker at the lower-right cell, the transient matrix of the other walker before the encounter here is also given as $M = P \cdot P^T$, which is symmetric.


\begin{figure}[hbp]
\centering
  \includegraphics[width=0.6\textwidth]{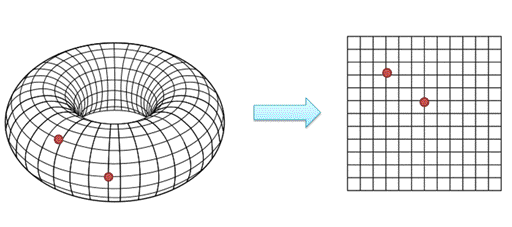}
  \caption{2-D Torus}
  \label{fig:3}
\end{figure}

Our main result is as follows: by suitably defining a Laplacian matrix $L$, the expected meeting time of the two walkers \textbf{E}$[\bm{\tau}]$ (i.e., the expectation of the first time that they meet on the same cell starting from the steady state uniform distribution) on a ring or torus could be explicitly expressed as the sum of the reciprocals of non-zero eigenvalues of $L$. We further conjecture based on empirical evidence that the result holds more generally for simple random walks (i.e., with equal transition probabilities) on arbitrary regular graphs.

\section{Method and Key Results}

\subsection{Preliminary}

Recall the standard definition of a Circulant Matrix: 

\begin{definition}[Circulant Matrix]
A circulant matrix is a matrix where each row vector is rotated one element to the right relative to the preceding row vector. A circulant matrix $A$ is fully specified by one vector, $\bm{a}$, which appears as the first row of $A$.
\end{definition}

\subsubsection{Properties of Circulant Matrix}
For arbitrary real, circulant matrix $A$ generated by $\{a_0,a_1,\cdots,a_{n-1}\}$with order $n$, we can find its eigenvalue in a general way following the approach indicated in~\cite{Kleinberg11}. First define vector $\bm{\xi_i}$ whose $j^{th}$ component is

\begin{equation} 
\xi_i(j) = \frac{1}{\sqrt{n}}w^{ij} \quad \text{where} \quad w=e^{\frac{2\pi}{n}} \text{is the $n^{th}$ roots of unity }
\end{equation}

We can prove the following properties:

(a) $<\xi_i,\xi_j>=\delta_{ij}$

(b) $\bm{A\xi_i} = \lambda_i \bm{\xi_i}, i = 1,2,...,n-1,0$\\

(a) shows that$\{\xi_i|i=1,2,\cdots,N\}$ are the orthogonal eigenvectors of $A$. $\lambda_i$ is the eigenvalue of $A$, which can be calculated by

\begin{eqnarray}
(A\xi_i)(j)& = &\sum_{k=1}^{n} A(j,k)\xi_i(k) \\
& = & \frac{a_0}{\sqrt{n}}w^{ij}+\frac{a_1}{\sqrt{n}}w^{i(j+1)}+\cdots+ \frac{a_{n-1}}{\sqrt{n}}w^{i(j+n-1)}\\
& = & \xi_i(j)(\sum_{k=0}^{n-1}a_k w^{ik})
\end{eqnarray}

Let $\lambda_i = \sum_{k=0}^{n-1}a_k w^{ik}$, then we have the property (b).

\begin{definition}[Block-Circulant Matrix]
If $A$ is a $n^2$-order partitioned circulant matrix generated by $A_0,A_1,\cdots,A_{n-1}$ where the $A_k$ are all $n$-order circulant matrices generated by $\{a_{i,0},a_{i,1},\cdots,a_{i,n-1}\}$ (see illustration below for a 9-order Block-Circulant Matrix). Then $A$ is called a block-circulant matrix. 
\end{definition}

$
 A =  \begin{bmatrix}
A_0 & A_1 & A_2 \\
A_2 & A_0 & A_1 \\
A_1 & A_2 & A_0
 \end{bmatrix}, 
 \quad \text{where} \quad
A_i = \begin{bmatrix}
a_{i,0} & a_{i,1} & a_{i,2} \\
a_{i,2} & a_{i,0} & a_{i,1} \\
a_{i,1} & a_{i,2} & a_{i,0}
 \end{bmatrix} \quad \text{for} \quad i=0,1,2.
$

\subsubsection{Properties of Block-Circulant Matrix}

Given index $i$, the coordinates of $i$ is $x_i=\text{quotient}(i-1,n),y_i=\text{remainder}（i-1,n）$。Then we need to modify the definition of $\bm{\xi_i}$ by

\begin{equation} 
\xi_x(j) = \frac{1}{n}w^{x_iy_i+x_jy_j} \quad \text{where} w=e^{\frac{2\pi}{n}} \text{is the $n^{th}$ roots of unity }
\end{equation} 

The properties given above in section 3.1.1 still hold, and we have \\ $\lambda_i = \sum_{l=0}^{n-1}\sum_{k=0}^{n-1} a_{l,k} w^{x_i l+y_i(k+1)}$ is the $i^{th}$ eigenvalue of $A$.

\subsection{Results on Circle}
\subsubsection{The Expected Meeting Time}
Let us first discuss the problem on the simplest graph, a 1-D circle.

\begin{theorem}
If two particles make independent random walks on a circle with an uniform initial distribution, then the expected meeting time is $\sum_{\lambda_i \ne 0} \lambda_i^{-1}$, where $\lambda_i$ is the $i^{th}$ eigenvalue of $L = I - PP^T$, and $P$ is the transition matrix for a single walker.
\end{theorem}

Put the transition probabilities in $M$ as the weight of edges. Then we get the Laplacian matrix,  

\begin{equation}       
L = I - M = I - PP^T
\end{equation}

which is a circulant matrix generated by $\{1-q_0,-q_1,-q_2,0,\cdots,0,-q_2,-q_1\}$, where $q_0 = p_1^2+p_2^2+p_3^2, q_1 = p_3(p_1+p_2), q_2 = p_1p_2$.

Let $T_i$,which is the $i^{th}$ component of vector $\bm{T}$, denote the expected meeting time with starting vertex $i$. Obviously, $T_0 = 0$. The initial distribution is $\bm{\pi}$. Then \textbf{E}$[\tau] = \bm{\pi}^T \bm{T}$. We can obtain a set of equations by recurrence:

\begin{equation}       
T_i=
  q_2 T_{i-2}+q_1T_{i-1}+q_0T_{i}+q_1T_{i+1}+q_2T_{i+2}+1 \quad i \ne 0
\end{equation}

Notice that the coefficients $q_0+2q_1+2q_2 = p_1^2+p_2^2+p_3^2 + 2p_3(p_1+p_2)+p_1p_2 = (p_1+p_2+p_3)^2=1$. By summing up the above equations, we have:

\begin{equation}       
T_0=
  q_2 T_{N-2}+q_1T_{N-1}+q_0T_{0}+q_1T_{1}+q_2T_{2}-(N-1)  
\end{equation}

Thus, the Laplacian matrix $L$ is the coefficient matrix of (4),(5). 

\begin{equation} 
L\bm{T}=\bm{\Delta t},\quad \text{where} \quad \Delta t = (1,1,1,\cdots,1,-(N-1))^T
\end{equation} 

Since L is a real, circulant matrix, we can use the conclusion in section 3.1.1. Taking the inner product of (9) with $\xi_i$ on both sides, from the symmetry of $L$ we have

\begin{equation} 
<\bm{LT},\bm{\xi_i}> = <\bm{T},\bm{L\xi_i}> = <\bm{T},\lambda_i\bm{\xi_i}> = \lambda_i (\sum_{k=1}^{N-1}T_k \frac{w^{ik}}{\sqrt{N}}+T_0)
\end{equation} 

\begin{equation} 
<\bm{\Delta t},\bm{\xi_i}> = \frac{1}{\sqrt{N}}( \sum_{k=1}^{N-1}w^{ik}-(N-1)) = 
\left\{
    \begin{array}{lr} 
    -\sqrt{N} & i\ne 0\\
    0 & i=0
    \end{array}
\right.
\end{equation} 

Notice that $\sum_{k=1}^{N-1}w^{ik} = -1$ for $i \ne 0$. Combined with (9), for $i \ne 0$,

\begin{equation} 
\sum_{k=1}^{N-1}\frac{T_k}{\sqrt{N}}w^{ik} = -\sqrt{N}(\lambda_i)^{-1}
\end{equation} 

Summing up by $i$, we have:

\begin{equation}
\sum_{i=1}^{N-1}\sum_{k=1}^{N-1}\frac{T_k}{\sqrt{N}}w^{ik}  = - \sqrt{N} \sum_{i=1}^{N-1} \lambda_i^{-1}
\end{equation}

\begin{equation}
\sum_{i=1}^{N-1}\sum_{k=1}^{N-1}T_k w^{ik} =  -N\sum_{i=1}^{N-1} \lambda_i^{-1}
\end{equation}

Changing the order of summation,

\begin{equation}
\sum_{k=1}^{N-1}T_k\sum_{i=1}^{N-1} w^{ik}  = -N\sum_{i=1}^{N-1} \lambda_i^{-1}
\end{equation}

\begin{equation}
\frac{1}{N}\sum_{k=1}^{N-1}T_k  = \sum_{i=1}^{N-1} \lambda_i^{-1}
\end{equation}

We assume the steady state distribution is the initial distribution. For any arbitrary regular graph, this is the uniform distribution. The expected meeting time is then given as:

\begin{equation}
\textbf{E}[\tau] = \bm{\pi}^T \bm{T} = \frac{1}{N}\sum_{k=1}^{N-1}T_k  = \sum_{i=1}^{N-1} \lambda_i^{-1}
\end{equation}

Note that this is the sum of the reciprocals of non-zero eigenvalues of $L$.
\subsubsection{The Order Estimation of $\textbf{E}[\tau]$}

For simplicity, we estimate the order of $\textbf{E}[\tau]$ for simple random walk (i.e., $p_1=p_2=p_3=\frac{1}{3}$):
\begin{equation}
\begin{split}
\textbf{E}[\tau]=&\sum_{i \ne 0}^N \frac{1}{3}(2-\frac{4}{3}\cos{\frac{\pi i}{N}}-\frac{2}{3}\cos{\frac{2 \pi i}{N}})^{-1}\\
=&\sum_{i=1}^N \frac{2}{9}(2-\cos{\frac{\pi i}{N}}-(\cos{\frac{\pi i}{N}})^2)^{-1}\\
=&\frac{2}{9}\sum_{i=1}^N \frac{1}{(2+t_i)(1-t_i)}\\
\end{split}
\end{equation}

where $t_i = \cos{\frac{\pi i}{N}} $. Thus $(2+t_i)^{-1} \in [1/3,1]$, which is bounded by constants. From~\cite{Montroll1969}, we have that summation $\sum_{i=1}^N (1-t_i)^{-1}$ is $O(N^2)$. Thus $\textbf{E}[\tau]$ is $O(N^2)$. On the other side, for $i = 1$, applying the Taylor Theorem we have

\begin{equation}
\frac{1}{1-t_1}=\frac{1}{1-\cos{\frac{\pi}{N}}} =\frac{1}{\Theta(1/N^2)} = \Theta(N^2)
\end{equation}

Thus $\textbf{E}[\tau]$ is also $\Omega(N^2)$, yielding that in fact for the 1-D circle, $\textbf{E}[\tau]$ grows with the size of the graph as $\Theta(N^2)$.

\subsection{Results on Torus}

\subsubsection{The Expected Meeting Time}
\begin{theorem}
If two particles make independent random walks on a torus with an uniform initial distribution, then the expected meeting time is $\sum_{\lambda_i \ne 0} \lambda_i^{-1}$, where $\lambda_i$ is the $i^{th}$ eigenvalue of $L = I - PP^T$, and $P$ is the transition matrix for a single walker.
\end{theorem}

Similarly, put the probabilities of transition in $M$ as the weight of edges. Then we get the Laplacian matrix.

\begin{equation}       
L = I - M =I - PP^T
\end{equation}

Let $T_i$ denotes the expected encounter time with starting point with index $i$, which is the $i^{th}$ component of vector $\bm{T}$. Obviously, $T_{N^2} = 0$. If the initial distribution is $\bm{\pi}$, then \textbf{E}$[\tau] = \bm{\pi}^T \bm{T}$. We can get a set of equations by recurrence (for a more readable notation here we write that $T_{x,y} = T_\text{Ind}(x,y)$). 

For ease of exposition, we illustrate below this recurrence equation for a simple random walk, that means the walker in the original model moves to its neighbour or stay still with the same probability $\frac{1}{5}$: 

\begin{equation}     
\begin{split}
T_{x,y}  = 
 & \frac{1}{25} T_{x\pm 2,y}+\frac{1}{25}T_{x,y\pm 2}+\frac{2}{25}T_{x\pm 1,y}+\frac{2}{25}T_{x,y\pm 1}\\
 +&\frac{2}{25}T_{x\pm 1,y\pm 1}+\frac{1}{5}T_{x,y}+1 \quad i \ne 0
  \end{split}
\end{equation}

Note that such a recurrence equation for $T_{x,y}$ could also be written for any random walk that moves to neighboring nodes with different probabilities. 

We also have:

\begin{equation} 
L\bm{T}=\bm{\Delta t},\quad \text{where} \quad \Delta t = (1,1,1,\cdots,1,-(N^2-1))^T
\end{equation} 

With the same approach in 3.2, we have

\begin{equation} 
<\bm{LT},\bm{\xi_i}> = <\bm{T},\bm{L\xi_i}> = <\bm{T},\lambda_i\bm{\xi_i}> = \lambda_i (\sum_{k=1}^{N} \sum_{l=1}^{N}\frac{T_{k,l}}{N}w^{x_i k+y_i l})
\end{equation} 

\begin{equation} 
<\bm{\Delta t},\bm{\xi_i}> = \frac{1}{N}( \sum_{(k,l)\ne (N,N)}w^{x_ik+y_il}-(N^2-1)) = 
\left\{
    \begin{array}{lr} 
    -N & i\ne 0\\
    0 & i=0
    \end{array}
\right.
\end{equation} 

Combined with (20) and $T_0=0$, summing up by $i$ for $i \ne 0$, we have

\begin{equation}
\sum_{i=1}^{N^2-1}\sum_{(k,l)\ne (N,N)}\frac{T_{k,l}}{N}w^{x_i k+y_i l}  = - N \sum_{i=1}^{N^2-1} \lambda_i^{-1}
\end{equation}

\begin{equation}
\sum_{(x,y)\ne (N,N)}\sum_{(k,l)\ne (N,N)}T_{k,l}w^{x_i k+y_i l}  = - N \sum_{i=1}^{N^2-1} \lambda_i^{-1}
\end{equation}

Change the sequence of summation, finally we have

\begin{equation}
\frac{1}{N^2}\sum_{(k,l)\ne (N,N)}T_{k,l}  = \sum_{i=1}^{N^2-1} \lambda_i^{-1}
\end{equation}

Note that we get actually the same expression as 1-D circle. Given the uniform initial distribution, the expected time $\textbf{E}[\tau]$ is the sum of the reciprocals of non-zero eigenvalues of $L$.

\subsubsection{The Order Estimation of $\textbf{E}[\tau]$}

Applied (6) to (25), we have

\[\textbf{E}[\tau]=
\sum_{\begin{subarray}{l} i,j=0 \\ (i,j)\ne(0,0)\end{subarray}}^{N-1}
\left( \frac{1}{25}
(20-2(\cos{\frac{4\pi i}{N}}+\cos{\frac{4\pi j}{N}})-4(\cos{\frac{2\pi i}{N}}+\cos{\frac{2\pi j}{N}})
-8\cos{\frac{2\pi i}{N}}\cos{\frac{2\pi j}{N}})
\right)^{-1}\]

which can be rewritten as
\begin{equation}
\textbf{E}[\tau]\equiv
\sum_{\begin{subarray}{l} i,j=0 \\ (i,j)\ne(0,0)\end{subarray}}^{N-1}
\frac{1}{2t_{ij}s_{ij}+3}\frac{1}{1-t_{ij}s_{ij}}
\end{equation}

where $t_{ij}=\cos{\frac{\pi(i+j)}{N}}, s_{ij}=\cos{\frac{\pi(i-j)}{N}}$. By applying the lemma(proved in \textbf{Appendix A}):

\begin{lemma}

If $\theta_1,\theta_2\in [0,\frac{\pi}{4}]$, then

\[\frac{1}{1-\cos{\theta_1}\cos{\theta_2}}\le \frac{4}{1-\cos{2\theta_1}\cos{2\theta_2}}\]

\end{lemma}

we can separate the summation into $\Theta(logN)$ parts, and prove that each part is $\Theta(N^2)$.
Thus finally we obtain that

\begin{equation}
\textbf{E}[\tau]~is~\Theta(N^2 log N)
\end{equation}

The complete proof is given in the \textbf{Appendix A}.

\section{Discussion}

\begin{figure}[!htb]
\centering
  \includegraphics[width=0.5\textwidth]{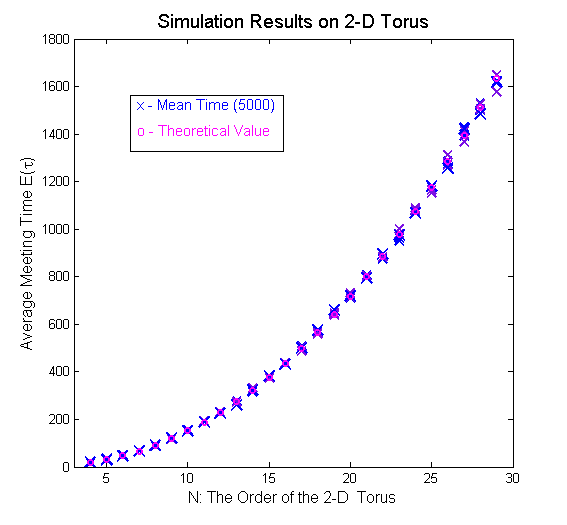}\\
  \caption{Simulation Results on 2-D Torus}
  \label{fig:4}
\end{figure}

We have proved that on the circle and the torus, the sum of the reciprocals of non-zero eigenvalues of $L = I - PP^T$ is the expected meeting time of two walkers. In fact, if the graph has a strong symmetry properties which guarantees $M = PP^T$ and $L$ is  (block-)circulant, then the proof still holds. The simulation results shown in figure~\ref{fig:4} match the conclusion in section 3.

\begin{figure}[!htb]
\centering
  \includegraphics[width=0.15\textwidth]{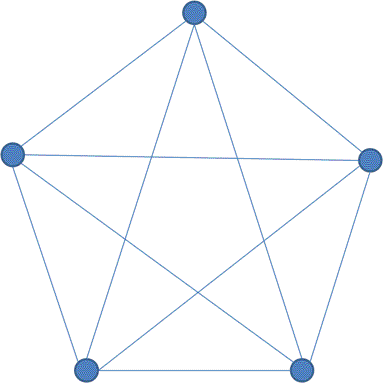} 
  \quad
  \includegraphics[width=0.17\textwidth]{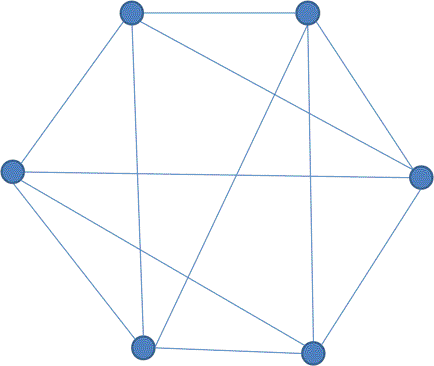} 
    \quad
  \includegraphics[width=0.17\textwidth]{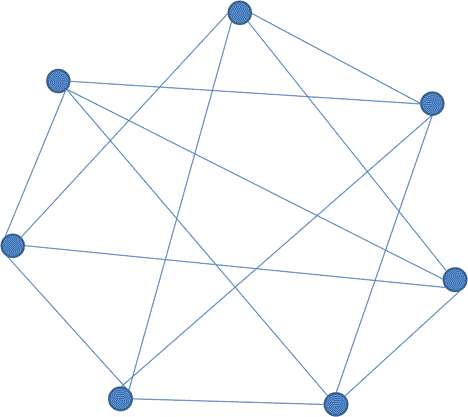}  
    \quad
  \includegraphics[width=0.15\textwidth]{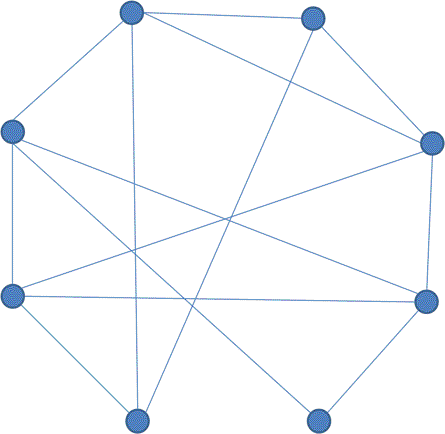} \\ 
  \caption{Special Cases for 4-regular Graph }
  \label{fig:5}
\end{figure}

Moreover, we find empirically that the expression even works for simple random walks on arbitrary regular graphs. This is not a trivial observation, since the symmetry of vertices doesn't hold for arbitrary regular graph, see the examples for 4-regular graphs in figure~\ref{fig:5}. In this case, the equivalent model approach of fixing one of the walkers at a particular location and defining the transition matrix of the other walker does not work. 

\begin{conjecture}[Expected Meeting Time on Regular Graph]

If two particles make independent \textbf{simple} random walks on a connected $d$-regular graph, and the initial distribution is uniform, then the expected meeting time $\textbf{E}[\tau]$ is $\sum_{\lambda_i \ne 0} \lambda_i^{-1}$, where $\lambda_i$ is the $i^{th}$ eigenvalue of $L = I - PP^T$, and $P$ is the transition matrix for a single walker.
\end{conjecture}

\begin{figure}[!htb]
\centering
  \includegraphics[width=0.45\textwidth]{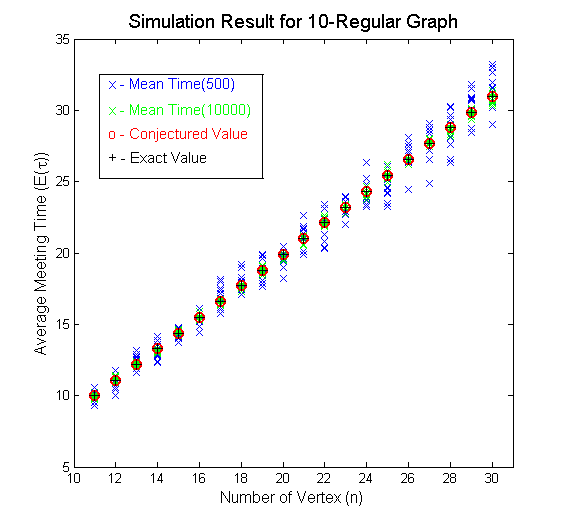} 
  \quad
  \includegraphics[width=0.45\textwidth]{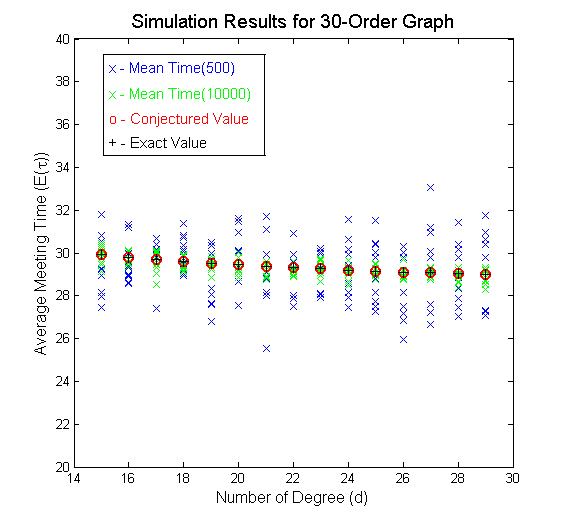} 
 \\ 
  \caption{Simulation Results on General Regular Graphs }
  \label{fig:6}
\end{figure}

Our conjecture is supported by empirical evidence which we present here. Figure~\ref{fig:6} shows  simulation results as well as relevant numerical calculations for simple random walks over arbitrary regular graphs. The left figure shows the results on 10-regular graphs, while the right one on graphs with 30 vertices. For each horizontal point, a single random graph is generated and fixed for averaging over multiple random initial conditions drawn from a uniform distribution. Each blue mark indicates the average meeting time when doing the experiment independently for 500 times, and green mark for 10000 times. The red mark indicates the conjectured value of the expected meeting time  (i.e. the sum of the reciprocals of non-zero eigenvalues of $L$). The black mark indicates the exact value of $\textbf{E}[\tau]$ which could be calculated by the definition of expectation once given transition probabilities (See Appendix B). In each case we see that the conjecture is valid.

One way to prove the conjecture may be to use the method in section 3; but for this approach we would need an additional conjecture.

\begin{conjecture}

If $A$ is the adjacency matrix of a connected d-regular graph $G$ with $n$ vertex, then $A$ has a set of orthogonal eigenvectors \{$\xi_1, \xi_2, \cdots, \xi_n$\} satisfying

$(a).\quad \xi_n = (1,1,\cdots, 1)^T$;

$(b).\quad \xi_i(n) = 1,$ for all $i$;

$(c). \quad  \sum_{j=1}^n \xi_i(j) = 0, $for all $i$;

$(d). \quad   \sum_{i=1}^n \xi_i(j) = 0, $for all $j$;

$(e). \quad   <\xi_i,\xi_j>=n \delta_{ij}$

\end{conjecture}

\begin{proposition}

\textbf{Conjecture 2} is a sufficient condition for \textbf{Conjecture 1}.

\end{proposition}

\textit{Proof.}
Suppose $\mu_1, ...,\mu_n$ is the eigenvalues of $A$.


We define a matrix $\tilde{L}$ as follows:
\begin{equation} 
\tilde{L} = I -  P \otimes P
\end{equation} 

where $P \otimes P$ is the kronecker product of P. Then from$P = (I+A)/(d+1)$, we have the eigenvalue of $P$ is $\beta_i = (\mu_i+1)/(d+1)$. Thus from the properties of kronecker product, the eigenvalue and eigenvector of $\tilde{L}$ is $\lambda_{i,j} =\beta_i \beta_j$ and $\xi_{i,j} = \xi_i \otimes \xi_j$.

We can similarly construct a recursive function of $T_{i,j}$, which indicates the expected meeting time with walkers on vertex $i$ and $j$. Obviously, $T_{i,i} = 0$. We can prove that $\tilde{L}T = \bm{\Delta t}$, where $\Delta t_{i,j} = 1$ if $i \ne j$, else $\Delta t_{i,i} = -(n-1)$. Then

\begin{equation} 
<\bm{\tilde{L} T},\bm{\xi_{i,j}}> = <\bm{T},\bm{L\xi_{i,j}}> = <\bm{T},\lambda_i\bm{\xi_{i,j}}> = \lambda_{i,j} \sum_{(k,l)=(1,1)}^{(n,n)}T_{k,l}\xi_i(k) \xi_j(l) 
\end{equation} 

Combined with (c) and (e) in \textbf{Conjecture 2}, we have

\begin{equation} 
\begin{split}
<\bm{\Delta t},\bm{\xi_{i,j}}> &=\sum_{k=1} ^n \left(-(n-1)\xi_i(k)\xi_j(k) + \xi_i(k)\sum_{l \ne k}\xi_j(l)\right)\\
& = \sum_{k=1} ^n -(n-1)\xi_i(k)\xi_j(k) - \xi_i(k)\xi_j(k)\\
& = -n <\xi_i,\xi_j> = -n^2 \delta_{ij}
\end{split}
\end{equation} 

Thus we have
\begin{equation}
\begin{split}
\sum_{(k,l)=(1,1)}^{(n,n)}T_{k,l}\xi_i(k) \xi_j(l) = \frac{1}{\lambda_{i,j}} <\bm{\Delta t},\bm{\xi_{i,j}}>= -n^2 \delta_{ij}
\end{split}
\end{equation}

Summing by $(i,j) \ne (n,n)$ and applying (d), finally we get the expression

\begin{equation}
\begin{split}
\textbf{E}[{\tau}] & = \frac{1}{n^2} \sum_{(i,j)=(1,1)}^{(n,n)}T_{i,j}
= \sum_{(i,j)=(1,1)}^{(n,n)} \delta_{ij} \frac{1}{\lambda_{i,j} } 
= \sum_i^n \frac{1}{\lambda_{i,i} }
\end{split}
\end{equation}

Notice that $\lambda_{i,i}$ is the same eigenvalue of $L = I - PP^T$ in our original definition of $L$. Thus we have proved that if \textbf{Conjecture 2} holds then the \textbf{Conjecture 1} would be true.

\begin{remark}

If we let $\xi_n$ be the eigenvector with eigenvalue $\mu=d$, then (a) holds.
\end{remark}

\begin{remark}
Since $\sum_{j=1}^n \xi_i(j) = (1,1,\cdots,1)^T \xi_i$, multiply $(1,1,\cdots,1)^T$ on the left of $L \xi_i = \lambda_i \xi_i$  and we have

\begin{equation}
\sum_{j=1}^n \xi_i(j) = \frac{1}{\lambda_i}\left((1,1,\cdots,1)^T L \right) \xi_i = 0
\end{equation}

Notice that the the row sum of $L$ is equal to 0. Thus we have (c).

\end{remark}

~\\
\textbf{\large Appendix A: The Proof for $\textbf{E}[\tau] = \Theta(N^2 log N)$ on 2-D Torus}

~\\

Recall \textbf{Lemma 1}:

\textit{
If $\theta_1,\theta_2\in [0,\frac{\pi}{4}]$, then
}
\[\frac{1}{1-\cos{\theta_1}\cos{\theta_2}}\le \frac{4}{1-\cos{2\theta_1}\cos{2\theta_2}}\]
~\\
\textit{Proof.}Let $s=\cos{\theta_1}$ and $t=\cos{\theta_2}$, then $\cos{2\theta_1}=2s^2-1$ and $\cos{2\theta_2}=2t^2-1$.

The inequality in lemma is equivalent to:

\begin{equation}
4-4st\ge 1-(2s^2-1)(2t^2-1)
\end{equation}

\[4ts-4t^2s^2+2t^2+2s^2\le 4\]

Let $f(t,s)=4ts-4t^2s^2+2t^2+2s^2$, since if $ts=c\le 1$ is fixed, $f$ attains its maximum at $t=1,s=c$. Thus, it remains to show $f(1,s)\le 4$, which is $-s^2+2s-1\le 0$, this inequality is correct and we complete the proof.
~\\

Recall the equation (26) which can be obtained by some trigonometric identities.

\[
\textbf{E}[\tau]=
\sum_{\begin{subarray}{l} i,j=0 \\ (i,j)\ne(0,0)\end{subarray}}^{N-1}
\frac{1}{2t_{ij}s_{ij}+3}\frac{1}{1-t_{ij}s_{ij}}
\]

Since $1\le 2ts+3\le 5$ for all $i,j$, then $\frac{1}{5}\le \frac{1}{2ts+3} \le 1$, which is bounded. Then we only need to estimate

\begin{equation}
\sum_{\begin{subarray}{l} i,j=0 \\ (i,j)\ne(0,0)\end{subarray}}^{N-1}\left(1-\cos{\frac{\pi(i+j)}{N}}\cos{\frac{\pi(i-j)}{N}}\right)^{-1}
\end{equation}

$(i,j)$ are uniformly distributed within the grid $[0,N]\times[0,N]$ (except the origin), then $(i\!+\!j,i\!-\!j)$ are uniformly distributed in a diamond area in $[0,2N]\times[-N,N]$, by the symmetry of cosine function and omitting a constant coefficient, it's equivalent to estimate

\begin{equation}
\sum_{\begin{subarray}{l}  p,q=0 \\ (p,q)\ne(0,0)\end{subarray}}^{N-1}\left(1-\cos{\frac{p\pi}{2N}}\cos{\frac{q\pi}{2N}}\right)^{-1}
\end{equation}

Since when we set $p = 0$(or $q = 0$), the summation is

\begin{equation}
\sum_{q=0}^{N-1}(1-\cos{\frac{q\pi}{2N}})^{-1}
\end{equation}

From~\cite{Montroll1969}, we have that summation is $O(N^2)$.

Thus, it remains to prove the following summation is in $\Theta(N^2 log N)$

\begin{equation}
\sum_{p,q=1}^{N-1}\left(1-\cos{\frac{p\pi}{2N}}\cos{\frac{q\pi}{2N}}\right)^{-1}
\end{equation}

Now let us partition the region into $\Theta(log N)$ parts, denote by

\begin{equation}
A_k = D_k \\ A_{k-1}, \quad \text{where} D_k ={(p,q)|1 \le p,q \le 2^k}, k = 0,1,2,\cdots, log N
\end{equation}

for all $k \ge 1$, $|A_k+1| = 4|A_k|$, and every term $(p,q)$ in $A_k$ corresponds to $(p,q),(p-1,q),(p,q-1),(p-1,q-1)$ in $A_k+1$.
Then applying the \textbf{Lemma 1} and the cosine function is non-negative and monotone decreasing in $[0,\pi/2]$, we can prove that 

\begin{equation}
S_k = \sum_{(p,q)\in A_k}\left(1-\cos{\frac{p\pi}{2N}}\cos{\frac{q\pi}{2N}}\right)^{-1}
\le \sum_{(p,q)\in A_{k+1}}\left(1-\cos{\frac{p\pi}{2N}}\cos{\frac{q\pi}{2N}}\right)^{-1} = S_{k+1}
\end{equation}

for $k = 0$, $S_0 \le S_1$ also holds by a simple calculation.

Notice that since $1- \cos{\frac{\pi}{2N}}=\Theta(\frac{1}{N^2})$, thus $S_1 =(1-\cos{\frac{\pi}{2N}}\cos{\frac{\pi}{2N}})^{-1}$ is $\Theta(N^2)$. The terms in $S_{logN}$ is bounded above by a constant $(1-\cos{\frac{\pi}{4}})^{-1} = 2$ and similarly bounded below by 0.5, then $S_{logN}$ is also $\Theta(N^2)$.

Thus, we have
\begin{equation}
\textbf{E}[\tau] = \sum_{k = 0}^{log N} S_k~is~\Theta(N^2 log N)
\end{equation}


~\\
\textbf{\large Appendix B: Calculating the Exact Value of $\bm{\textbf{E}[\tau]}$}

~\\
The exact value of expected meeting time could be calculated in the following way:
Suppose there are two walkers $a$ and $b$. We denote the state that $a$ is at vertex $i$ while $b$ is at vertex $j$ by $S_{(i,j)}$, with index $((i-1)N+j)$. Thus if the transition matrix for a single walker is $P$, then the transition matrix for the states of two walkers is $Q = P \otimes P$ except for the $(i-1)N+i^{th}$ rows(the absorbing states),  which are all zeros expect the the $iN+i^{th}$ component. Let $\Lambda = \{(i-1)N+i|i=1,2,\cdots,N\}$, and $S_\Lambda$ is the set of absorbing states. $S(\tau)$ indicates the state at time $\tau$.

Recall the definition of expectation, we have
\begin{equation}
\textbf{E}[\tau]  =\sum_{\tau =0}^{\infty} \tau Pr[S(\tau) \in S_\Lambda,S(\tau - 1) \notin S_\Lambda]
\end{equation}

that equals to

\begin{equation}
\begin{split}
\textbf{E}[\tau] =& \sum_{\tau =0}^{\infty} \tau \sum_{k \in S_\Lambda} \sum_{l \notin S_\Lambda } Pr[S(\tau) =k,S(\tau - 1) =l]\\
=& \sum_{\tau =0}^{\infty} \tau \sum_{k \in S_\Lambda} \sum_{l \notin S_\Lambda } Pr[S(\tau) =k|S(\tau - 1)=l]Pr[S(\tau - 1) =l]\\
=& \sum_{\tau =0}^{\infty} \tau \sum_{k \in S_\Lambda} \sum_{l \notin S_\Lambda } Pr[S(\tau) =k|S(\tau - 1) =l]\bm{p_0}\cdot Q^{\tau-1} \cdot \bm{e_l}\\
=& \bm{p_0} \sum_{\tau =0}^{\infty}(\tau \cdot Q^{\tau-1} )\cdot \bm{b}
\end{split}
\end{equation}

where $\bm{b}$ is a column vector with $n^2$ component, $\bm{b}(j) = \sum_{i \in \Lambda}Q(i,j)$ if $j \notin \Lambda$, $\bm{b}(j) = 0$ if $j \in \Lambda$. Then applying the series summation approach to matrix, finally we have

\begin{equation}
\textbf{E}[\tau]  = \bm{p_0}(I-B)^{-2}\bm{\tilde{b}}.
\end{equation}

where $B$ is the sub-matrix of $Q$ deliminating the rows and columns with index in $\Lambda$, $\tilde{b}$ is the sub-vector of $b$ deliminating the rows and columns with index in $\Lambda$.



\end{document}